\documentclass{article}

\newtheorem{thm}{Theorem}[section]
\newtheorem{cor}[thm]{Corollary}

\newtheorem{lem}[thm]{Lemma}

\newcommand{\be}{\begin{equation}}
\newcommand{\ee}{\end{equation}}
\newcommand{\ben}{\begin{enumerate}}
\newcommand{\een}{\end{enumerate}}
\newcommand{\pa}{{\partial}}
\newcommand{\R}{{\rm R}}

\newcommand{\g}{{\bf g}}
\newcommand{\pxi}{{\pa \over \pa x^i}}
\newcommand{\pxj}{{\pa \over \pa x^j}}

\newcommand{\qed}{\hspace*{\fill}Q.E.D.}  
\title{\Large On Negatively Curved  Finsler Manifolds\\
 of Scalar Curvature}
\author{Xiaohuan Mo\footnote{supported by the National Natural Science Foundation of China (10171002)} and Zhongmin Shen}

\date{February 15, 2003}

\begin{document}

\maketitle

\begin{abstract}
In this paper, we prove a global rigidity theorem for negatively curved Finsler metrics on a compact manifold of dimension $n \geq 3$.
 We show that for such a Finsler manifold, if the flag curvature is a scalar function on the tangent bundle, then  the Finsler metric is of Randers type. We also study the case when the Finsler metric is locally projectively flat. \end{abstract}

\bigskip

\section{Introduction}
The flag curvature in Finsler geometry is a natural extension of the sectional curvature in Riemannian geometry, which is first introduced by L. Berwald \cite{Be1}\cite{Be2}. 
 For a Finsler manifold $(M, F)$, the flag curvature is  a function ${\bf K}(P, y)$ of tangent planes $P\subset T_xM$ and  directions $y\in P$. 
$F$  is said to be 
 {\it of scalar curvature} if 
the flag curvature ${\bf K}(P, y)={\bf K}(x, y)$ is independent of 
flags $P$ associated with any fixed flagpole $y$. 
Finsler metrics of scalar curvature are the natural extension of Riemannian metrics of isotropic sectional curvature (which are of  constant sectional curvature in dimension $n\geq 3$ by the Schur Lemma).
One of the important problems in Finsler geometry is to characterize
 Finsler manifolds of scalar curvature. 
It is known that locally projectively flat Finsler metrics  are of scalar curvature  \cite{Be3}\cite{Be4}. The converse is false.  There are lots of  non-projectively flat Finsler metrics  of scalar curvature.  To our best knowledge,  all  known non-projectively flat  Finsler metrics of scalar curvature are in the form $F = \alpha +\beta$, where 
$\alpha$ is a Riemannian metric  and $\beta$ is a $1$-form \cite{BaRo}\cite{Sh2}\cite{Sh3}.  Such metrics are called {\it Randers metrics}. Randers metrics arise naturally from physical applications \cite{Ra}. 
 In this paper, we prove the following global result.

\begin{thm}\label{thmC-reducible} Let $(M, F)$ be a compact negatively curved Finsler manifold of dimension $n\geq 3$. Suppose that 
$F$ is of  scalar curvature.  
Then $F$ is a  Randers metric. 
\end{thm}

Theorem \ref{thmC-reducible} greatly narrows down the possibility of compact negatively curved Finsler manifolds of scalar curvature. 
Although Randers metrics are very special Finsler metrics, the classification of  Randers metrics of scalar curvature has not been completely done yet. See \cite{YaSh} \cite{Ma2} \cite{MaSh}  \cite{BaRo} \cite{BaRoSh} for the classification of 
Randers metrics of constant flag curvature.

It is known  that  the Funk metric and the Hilbert metric  on a strongly convex domain in $\R^n$ 
have  negative constant flag curvature, but they are not Randers metrics when the domain is not an ellipsoid \cite{Be3}\cite{Be4}\cite{Fk}. Thus the compactness of the manifold in Theorem \ref{thmC-reducible} can not be replaced by completeness of the metric without additional growth condition on some non-Riemannian quantities.

 A Randers metric is locally projectively flat if and only if $\alpha$ is of constant sectional curvature and $\beta$ is a closed $1$-form. This is a direct consequence of a result in \cite{BaMa} and the Beltrami theorem on projectively flat Riemannian metrics. 
Thus if  $(M, \alpha)$ is a compact hyperbolic Riemannian manifold 
and $\beta$ is an arbitrary closed $1$-form on $M$, then 
the Randers metric $F = \alpha + \beta$ is of scalar curvature  with negative flag curvature, provided that the $\alpha$-norm of $\beta$ is sufficiently small. If in addition, $F$ has constant flag curvature, then $\beta=0$. This is because that $F$ is projectively equivalent to $\alpha$ and the Riemannian universal cover $(\tilde{M}, \tilde{\alpha})$ is isometric to the Klein model (\cite{Be4}\cite{Fk}\cite{Sh1}). This fact also follows from 
Akbar-Zadeh's rigidity theorem: every Finsler metric on a compact manifold of negative constant flag curvature must be Riemannian \cite{AZ}. 
Since locally projecively flat Finsler metrics are of scalar curvature, based on Theorem \ref{thmC-reducible} and the above arguments,   we obtain the following

\begin{thm} Let $(M, F)$ be a  compact negatively curved  Finsler manifold of dimension $n\geq 3$. $F$ is locally projectively flat if and only if $F = \alpha +\beta $ is a Randers metric where 
$\alpha$ is of constant sectional curvature and $\beta$ is a closed $1$-form.
\end{thm}

Besides the flag curvature, there are several important non-Riemannian curvatures. 
Let $(M, F)$ be a Finsler manifold. The second derivatives of ${1\over 2} F_x^2$ at $y\in T_xM\setminus\{0\}$ is an inner product $\g_y$ on $T_xM$.
 The third order derivatives of ${1\over 2} F_x^2$ at  $y\in T_xM\setminus\{0\}$ is a symmetric trilinear forms ${\bf C}_y$ on $T_xM$. We call $\g_y$ and ${\bf C}_y$ the {\it fundamental form} and  the {\it Cartan torsion}, respectively. 
The rate of change of the Cartan torsion along geodesics is the {\it Landsberg curvature} (or the {\it $L$-curvature}) $ {\bf L}_y$ on $T_xM$
for any $y\in T_xM\setminus\{0\}$.
Set ${\bf J}_y:= \sum_{i=1}^n {\bf L}_y(e_i, e_i, \cdot )$,
where $\{e_i\}$ is an orthonormal basis for $(T_xM, \g_y)$. 
${\bf J}_y$ is called the {\it mean Landsberg curvature} (or the {\it J-curvature}). 
$F$ is said to be {\it  Landsbergian} if ${\bf L}=0$,  and {\it weakly Landsbergian} if ${\bf J}=0$. See \cite{Sh4}\cite{Sh5}.

\begin{cor}\label{cor1}
Let $(M, F)$ be a compact  negatively curved  Finsler  manifold of dimension $n\geq 3$. Suppose that $F$ is  of scalar curvature and weakly Landsbergian. Then 
$F$ is Riemannian.
\end{cor}

According to Numata \cite{Nu}, any   Landsberg metric of scalar curvature with non-zero flag curvature is Riemannian. Here we weaken Numata's condition on 
the Landsberg curvature and impose the compactness on the manifold instead. 
We do not know if Corollary \ref{cor1} is still true  for  weakly Landsberg metrics of scalar curvature with  positive flag curvature. 
It is known that every weakly Landsberg metric of non-zero {\it constant} flag curvature is Riemannian (see Theorem 9.1.1  in \cite{Sh4}). 

\bigskip

There is another important quantity defined using the spray of $F$.
We call ${\bf E}_y$
the {\it mean Berwald curvature} (or the  {\it E-curvature}). 
$F$ is said to be {\it weakly Berwaldian} if ${\bf E}=0$.  See \cite{Sh4}\cite{Sh5}.

\begin{cor}\label{cor2}
Let $(M, F)$ be a compact negatively curved  Finsler manifold of dimension $n\geq 3$. Suppose that $F$ is of scalar curvature and weakly Berwaldian.
Then $F$ is Riemannian.
\end{cor}

The compactness in Corollary \ref{cor2} can not be dropped. Example 3.3.3 in \cite{BaRo} is a weakly Berwaldian  Randers metric with negative constant curvature. 

\section{Preliminaries}

In this section, we are going to 
give a brief description on several quantities in Finsler geometry.

Let $(M, F)$ be a Finsler manifold of dimension $n$. Fix a local frame 
$\{ {\bf b}_i\}$ for $TM$. The Finsler metric  $F = F( y^i {\bf b}_i )$ is a 
function of $(x^i, y^i)$. Let 
\[ g_{ij}(x,y): = {1\over 2} [F^2]_{y^iy^j}(x,y), \ \ \ \ \ C_{ijk}(x, y):= {1\over 4} [F^2]_{y^iy^jy^k}(x, y).\] 
For a non-zero vector $y=y^i {\bf b}_i\in T_xM$, the fundamental  form $\g_y$ on $T_xM$ is a bilinear symmetric form defined by $\g_y({\bf b}_i, {\bf b}_j):= g_{ij}(x,y)$, and the Cartan torsion ${\bf C}_y$  on $T_xM$ is a trilinear symmetric form on $T_xM$ defined by ${\bf C}_y ({\bf b}_i, {\bf b}_j, {\bf b}_k ) := C_{ijk}(x, y)$. The {\it mean Cartan torsion} ${\bf I}_y$ is  a linear form on $T_xM$ defined by 
\[{\bf I}_y({\bf b}_i) = I_i (x,y):= g^{jk}(x,y) C_{ijk}(x, y).\]
Let 
\[
M_{ijk}:=C_{ijk} - {1\over n+1} \Big \{ I_i h_{jk} + I_j h_{ik} + I_k h_{ij} \Big \},\label{Matsumoto}
\]
where 
\[ h_{ij} := FF_{y^iy^j} = g_{ij}- {1\over F^2} g_{ip}y^p g_{jq} y^q .\]
We obtain a symmetric  trilinear form 
${\bf M}_y $ on $T_xM$ defined by 
${\bf M}_y({\bf b}_i, {\bf b}_j, {\bf b}_k) := M_{ijk}(x, y)$. 
This quantity is introduced by M. Matsumoto \cite{Ma1}. Thus we call ${\bf M}_y$ the {\it Matsumoto torsion}. 
Matsumoto proves that every Randers metric satisfies that ${\bf M}_y=0$.
Later on, Matsumoto-H\={o}j\={o} proves that the converse is true too.
\begin{lem}\label{lemMaHo}{\rm (\cite{Ma1}\cite{MaHo})}
A Finsler metric $F$ on a manifold of dimension $n\geq 3$ is a Randers metric if and only if ${\bf M}_y =0$, $\forall y\in TM\setminus\{0\}$.
\end{lem}

Finsler metrics in this paper are always assumed to be regular in all directions. If  this regularity is not imposed, Matsumoto-H\={o}j\={o}'s theorem says that $F$ has vanishing Matsumoto torsion if and only if $F=\alpha+\beta$ or $F= \alpha^2/\beta$, where $\alpha$ is a Riemannian metric and $\beta$ is a $1$-form on $M$.

\bigskip
The spray  of $F$ is a vector field on $TM\setminus\{0\}$. In a standard local coordinate system $(x^i, y^i)$ in $TM$, the spray is given by
\[ {\bf G} = y^i\pxi - 2 G^i(x, y) {\pa \over \pa y^i},
\]
where $G^i(x, y):= {1\over 4} g^{il}(x, y) \{ [F^2]_{x^ky^l}y^k-[F]^2_{x^l}   \}$.
A Finsler metric $F$ is called a {\it Berwald metric} if $G^i(x, y)
= {1\over 2} \Gamma^i_{jk}(x)y^jy^k$ are quadratic in $y\in T_xM$. 
For a  vector $y\in T_xM\setminus\{0\}$, define 
a symmetric bilinear form ${\bf E}_y$ on $T_xM$ by
\[ {\bf E}_y \Big ( \pxi|_x, \pxj|_x \Big )= E_{ij}(x, y) := {1\over 2} {\pa^2 \over \pa y^i\pa y^j} \Big [ {\pa G^m \over \pa y^m} \Big ].\]
We call ${\bf E}_y$ the {\it mean Berwald curvature} (or the {\it E-curvature}). Clearly, for any Berwald metric, ${\bf E}_y=0$. In general, the converse is not true. Thus Finsler metrics with vanishing E-curvature are called weakly Berwald metrics. 
See \cite{Sh4}.

Let $c(t)$ be a $C^{\infty}$ curve and $U(t)= U^i(t) \pxi|_{c(t)}$ be a vector field along $c$. Define the covariant derivative of $U(t)$ along $c$ by
\[ D_{\dot{c}} U(t) := \Big \{ {dU^i\over dt}(t) + U^j (t) {\pa G^i\over \pa y^j} \Big ( c(t), \dot{c}(t) \Big ) \Big \} \pxi_{|c(t)}.\]
$U(t)$ is said to be {\it linearly parallel} if $D_{\dot{c}} U(t)=0$.

For a vector $y\in T_xM$, define 
\begin{eqnarray*}
{\bf L}_y (u,v,w): & = & {d\over dt} \Big [ {\bf C}_{\dot{\sigma}(t) } \Big ( U(t), V(t), W(t) \Big )\Big ]|_{t=0},\\
{\bf J}_y (u): & = & {d\over dt} \Big [ {\bf I}_{\dot{\sigma}(t) } \Big ( U(t) \Big )\Big ]|_{t=0},
\end{eqnarray*}
where   $\sigma(t)$ is the geodesic with $\sigma(0)=x$, $\dot{\sigma}(0)=y$
and $U(t), V(t), W(t)$ are  linearly parallel vector fields along $\sigma$ with 
$U(0)=u, V(0)=v, W(0)=w$. 
We call ${\bf L}_y$ the {\it Landsberg curvature} (or the {\it L-curvature}). The Landsberg curvature measures the rate of change of the Cartan torsion along geodesics.
Let  $L_{ijk}(x,y):= {\bf L}_y({\bf b}_i, {\bf b}_j, {\bf b}_k)$ and $J_i(x, y):= {\bf J}_y({\bf b}_i)$. We have that $J_i (x, y)= g^{jk}(x, y) L_{ijk}(x, y)$. Thus we call
${\bf J}_y$  the {\it mean Landsberg curvature} (or the {\it J-curvature}). See \cite{Sh4}.

The above mentioned  quantities all vanish for Riemannian metrics, hence they do not appear in Riemannian geometry. The Riemann curvature plays the most important role  in Finsler or Riemannian geometry. The Riemann curvature is introduced by Riemann in 1854 for Riemannian metrics and extended to Finsler metrics by L. Berwald in 1926 \cite{Be1}\cite{Be2}. A quick definition using Riemannian geometry is given as follows. For a non-zero vector $y\in T_xM$, take a local extension $Y$ of $y$ such
all integral curves of $Y$ are geodesics of $F$. Then $h:=\g_Y$ is a Riemannian metric. Let ${\rm R}(u, v)w$ denote the Riemannian curvature tensor of $h$. The {\it Riemann curvature} is a linear map  ${\bf K}_y: T_xM \to T_xM$ defined by 
\[
{\bf K}_y(u):= {\rm R}(u, y)y, \ \ \ \ \  u\in T_xM.\]
The Riemann curvature satisfies
\[ {\bf K}_y(y)=0, \ \ \ \ \g_y({\bf K}_y (u), v) = \g_y(u, {\bf K}_y(v)).\]
 For a  tangent plane $P\subset T_xM$ and a unit vector $y\in P$, the following number \[
{\bf K}(P, y):= \g_y({\bf K}_y (y^{\bot}), y^{\bot} )\]
 is called the {\it flag curvature} of the ``flag'' $(P, y)$, where 
$y^{\bot}\in P$ is an unit vector such that $\g_y(y, y^{\bot})=0$. See \cite{Sh4}\cite{Sh5}.

\bigskip
To find the relationship among the above quantities, we need the Bianchi identities.
Lifting  the local frame $\{ {\bf b}_i \}$ to a local frame $\{ {\bf e}_i\}$ for $\pi^*TM$ by setting ${\bf e}_i(x,y):= (y, {\bf b}_i(x))$. 
Let $\{\omega^i, \omega^{n+i}\}$ denote the corresponding local coframe for $T^*(TM\setminus\{0\})$. The above mentioned quantities can be viewed as tensors on $TM\setminus\{0\}$. For example, the Cartan torsion can be expressed as
${\bf C}= C_{ijk} (x, y) \omega^i \otimes \omega^j \otimes \omega^k$
and the Riemann curvature can be expressed as ${\bf K}= K^i_{\ k}(x, y) \omega^k \otimes {\bf e}_i$. 

The Chern connection forms are the unique local $1$-forms $\omega_j^{\ i} $ satisfying 
\[ d\omega^i = \omega^j \wedge \omega_j^{\ i},\]
\[dg_{ij} = g_{ik} \omega_j^{\ k} + g_{kj} \omega_i^{\ k} + 2 C_{ijk} \omega^{n+k},
\]
\[\omega^{n+k} = dy^k + y^j \omega_j^{\ k},\]
where $y^i$ are viewed as local functions on $TM$ whose values $y^i$ at $y$ are defined by $y= y^i {\bf b}_i$. See \cite{Ch}\cite{BCS}. With the Chern connection, we define covariant derivatives of quantities on $TM$  in the usual way. For example, for a scalar function $f$, we define $f_{|i}$ and $ f_{\cdot i}$ by
\[ df = f_{|i} \omega^i + f_{\cdot i} \omega^{n+i},\]
for the mean Cartan torsion ${\bf I}=I_i \omega^i$, define $I_{i|j}$ and $I_{i\cdot j}$ by
\[ d I_i - I_k \omega_i^{\ k} = I_{i|j}\omega^j + I_{i\cdot j} \omega^{n+j}.\]

For a tensor ${\bf T}= T_{i \cdots k } \omega^i \otimes \cdots \otimes \omega^k$, we have
\[  T_{i \cdots k \; \cdot m} = {\pa T_{i\cdots k}\over \pa y^m}.\]
For a non-zero vector $y\in T_xM$, the tensor ${\bf T}$ induces  a multi-linear form  ${\bf T}_y (u, \cdots, w):= 
T_{i \cdots k} (x, y) u^i \cdots w^k$ on $T_xM$. 
Let $\sigma(t)$ denote the geodesic with $ \dot{\sigma}(0)=y$. We have 
\[ {d \over dt} \Big [ {\bf T}_{\dot{\sigma}(t) } \Big ( U(t), \cdots,  W(t) \Big )\Big ]
=  T_{i\cdots k |m} (\sigma(t), \dot{\sigma}(t) )\dot{\sigma}^m (t)U^i(t) \cdots W^k(t),\]
where $ U(t)=U^i(t) \pxi|_{\sigma(t)} , \cdots, W(t)=W^k(t) {\pa \over \pa x^k}|_{\sigma(t)} $ are linearly parallel vector fields along $\sigma$. 
Thus the L-curvature ${\bf L}= L_{ijk} \omega^i\otimes \omega^j \otimes \omega^k$ and the J-curvature ${\bf J}= J_i \omega^i$ are given by
\be
L_{ijk} = C_{ijk|m}y^m, \ \ \ \ \ \ J_i = I_{i|m}y^m.\label{LJE}
\ee
Express the Riemann curvature by
${\bf K}= K^i_{\ k} \omega^k \otimes {\bf e}_i$.
Let 
\be
 K^i_{\ kl}:= {1\over 3} \Big \{ K^i_{\ k \cdot l} - K^i_{\ l\cdot k} \Big \}.\label{Kikl}
\ee
In fact,  $K^i_{\ kl}$ and $L^i_{\ kl}:= g^{ij}L_{jkl}$ are determined by the following equation
\[ \Omega^i := d\omega^{n+i}-\omega^{n+j} \wedge \omega_j^{\ i}
= {1\over 2} K^i_{\  kl}\omega^k \wedge \omega^l - L^i_{\ kl} \omega^k \wedge \omega^{n+l}.\]
But we do not need this fact here. 
According to (10.13) in \cite{Sh5}, the Landsberg curvature satisfies the following identities
\be
L_{ijk;m}y^m + C_{ijm}K^m_{\ \ k} 
= - {1\over 2} g_{im} K^m_{\ \ kl\cdot j} y^l - {1\over 2} g_{jm} K^m_{\ \ kl\cdot j} y^l,\label{LCKK}
\ee
Here $L_{ijk;m}$ denotes the horizontal covariant derivative with respect to the Berwald connection and $K^i_{\ \ kl\cdot j}= R^{\ i}_{j \ kl}$ are the coefficients of the Riemannian curvature (or the hah-curvature) of the Berwald connection.
 Since the Chern connection and the Berwald connection differ by the Landsberg curvature, we have $L_{ijk|m}y^m = L_{ijk;m}y^m$.
Plugging (\ref{Kikl}) into (\ref{LCKK}) yields
\begin{eqnarray}
L_{ijk|m}y^m + C_{ijm}K^m_{\ \ k} & = &  - {1\over 3}g_{im}K^m_{\ \ k\cdot j}
- {1\over 3} g_{jm} K^m_{\ \ k\cdot i}\nonumber\\
& &  - {1\over 6} g_{im}K^m_{\ \ j\cdot k}
- {1\over 6} g_{jm}K^m_{\ \ i\cdot k}. \label{Moeq1}
\end{eqnarray} 
Contracting (\ref{Moeq1}) with $g^{ij}$
gives
\be
J_{k|m}y^m + I_mK^m_{\ \ k}  = -  {1\over 3}\Big \{ 2 K^m_{\ \  k\cdot m} +  K^m_{\ \  m\cdot k}\Big \} . \label{Moeq2}
\ee
Equations (\ref{Moeq1}) and (\ref{Moeq2}) are established in \cite{Mo1} in a different way.

\section{Proof of Theorem \ref{thmC-reducible}}

In this section, we are going to prove a generalization of Theorem \ref{thmC-reducible}. 
First we define the norm of the Matsumoto torsion  at $x\in M$ by \index{norm of the Matsumoto torsion}
\[
 \| {\bf M}\|_x: = \sup_{y, u, v, w\in T_xM\setminus\{0\}} { F(y)|{\bf M}_y (u, v, w) | \over \sqrt{\g_y (u, u) \g_y(v,v)\g_y(w,w)} }.\label{normM}
\]
The Matsumoto torsion  is said to {\it  grow sub-exponentially} at rate of $k >0$ if 
there is any point $x\in M$ such that  as $r \to \infty$, 
\[ M(x, r):= \sup_{\min ( d(z, x), d(x, z) ) \leq r} \| {\bf M}\|_z  
= o ( e^{kr} ).    
\]

\begin{thm} \label{thmC-reducible2} Let $(M, F)$ be an $n$-dimensional complete  Finsler manifold of scalar curvature with flag curvature ${\bf K}\leq -1$ ($n\geq 3$). Suppose that the Matsumoto torsion grows sub-exponentially at rate of $k=1$. Then $F$ is a Randers metric. In particular, if $(M, F)$ is an $n$-dimensional compact  Finsler manifold of scalar curvature with flag curvature ${\bf K}\leq -1$, then it is a Randers metric.
\end{thm}
{\it Proof}:  We will first prove that the Matsumoto torsion vanishes. To prove this, we assume that the Matsumoto torsion ${\bf M}_y (u, u, u)= M_{ijk}(x, y) u^iu^ju^k \not=0$ for some $y, u\in T_xM \setminus\{0\}$ with $F(x, y)=1$. 
Let $\sigma(t)$ be the unit speed geodesic with $\sigma(0)=x$ and $\dot{\sigma}(0)=y$. 
Let $U(t)$ denote the linear parallel vector field along $\sigma$, that is,
$D_{\dot{\sigma}}U(t) =0$.
From the above equation, we see that a linearly parallel vector field 
$U(t)$ along $\sigma$ linearly depends on its initial value $U(t_o)$ at a point $\sigma(t_o)$. 

Let 
\[ {\cal M}(t):= {\bf M}_{\dot{\sigma}(t)} \Big ( U(t), U(t), U(t) \Big )
= M_{ijk} \Big ( \sigma(t), \dot{\sigma}(t)\Big ) U^i(t) U^j(t) U^k(t).\]
We have 
\[ {\cal M}''(t)=M_{ijk|p|q}\dot{\sigma}^p(t) \dot{\sigma}^q (t) \Big ( \sigma(t), \dot{\sigma}(t)\Big ) U^i(t) U^j(t) U^k(t).\]

Now we assume that $F$ is of scalar curvature  with flag curvature ${\bf K}= {\bf K}(x, y)$. This is equivalent to the following identity:
\be
K^i_{\ k} = {\bf K} F^2 \; h^i_k, \label{Kikiso1}
\ee
where $h^i_k := g^{ij} h_{jk}$. 
Differentiating (\ref{Kikiso1}) yields
\be
 K^i_{\ k\cdot l}
=  {\bf K}_{\cdot l} F^2 \; h^i_k + {\bf K} \Big \{ 2 g_{lp}y^p \delta^i_k - g_{kp}y^p  \delta^i_l- g_{kl} y^i  \Big \}.\label{MKdiff}\ee
By  (\ref{Moeq1}),  (\ref{Moeq2}) and (\ref{MKdiff}), we obtain 
\be
 L_{ijk|m}y^m =  - {1\over 3}F^2 \Big \{ {\bf K}_{\cdot i} h_{jk}   + {\bf K}_{\cdot j} h_{ik}  + {\bf K}_{\cdot k} h_{ij} + 3 {\bf K} C_{ijk}\Big \} \label{AZeq1}
\ee
and
\be
J_{k|m}y^m = - {1\over 3}F^2\Big \{ (n+1) {\bf K}_{\cdot k} + 3{\bf K} I_k \Big \}.\label{AZeq2}
\ee
By (\ref{LJE}), we have 
\[ C_{ijk|p|q}y^py^q = L_{ijk|m}y^m, \ \ \ \ \ I_{k|p|q}
y^py^q = J_{k|m}y^m.\]
\be
M_{ijk|p|q} = L_{ijk|m}y^m - {1\over n+1} \Big \{ J_{i|m}y^m h_{jk}
+ J_{j|m}y^m h_{ik} + J_{k|m}y^m h_{ij} \Big \}.
\label{AZeq3}
\ee
Plugging  (\ref{AZeq1}) and (\ref{AZeq2}) into (\ref{AZeq3}) yields
\be
 M_{ijk|p|q} y^p y^q + {\bf K} F^2 M_{ijk} =0. \label{Sijk}
\ee

It follows from (\ref{Sijk}) that 
\[
{\cal M}''(t) + {\bf K}(t) {\cal M}(t)=0,
\]
where ${\bf K}(t):= {\bf K}\Big (\sigma(t), \dot{\sigma}(t) \Big ) \leq -1$.
By an elementary argument, we have
\be
 \Big | {\cal M}(t) \Big | \geq \Big | {\cal M}(0)\cosh (t)+ {\cal M}'(0) \sinh(t) \Big |, \ \ \ \ \ \alpha < t < \beta,\label{Mvc}
\ee
where $(\alpha, \beta)$ is the maximal interval containing $0$ such that 
 the function on the right side of (\ref{Mvc}) is not equal to zero.
Let 
\[ M(x, t):= \sup_{ \min\{ d(z, x), d(x, z) \} \leq t}
\|{\bf M}\|_{z} .\]
Since 
$d ( \sigma(-t), x) \leq t$ and $ d (x, \sigma(t) ) \leq t$ for any $t >0$, 
we have
\[ M(x, r) \geq \max\Big \{ {\cal M}(t) \ \Big | \ |t| \leq r \Big \}.\]

Suppose that ${\cal M}'(0)=0$ or has the same sign as ${\cal M}(0)$. 
Then $\beta = \infty$ and 
\[ M(x, r) \geq \Big | {\cal M}(r) \Big | \geq  \Big | {\cal M}(0) \Big | \cosh (r)+ \Big | {\cal M}'(0)\Big |  \sinh(r), \ \ \ \ \ \ r >0.\]
Suppose that ${\cal M}'(0)$ has the opposite sign as ${\cal M}(0)$. 
Then $\alpha =-\infty$ and 
\[ M(x, r) \geq \Big | {\cal M}(-r) \Big | 
\geq \Big | {\cal M}(0)\Big  | \cosh(r) + \Big | {\cal M}'(0)\Big |  \sinh(r), \ \ \ \ \ \ r >0.\]
In either case, 
\[ \liminf_{r\to \infty} { M(x, r) \over e^{r} }\geq |{\cal M}(0) | + |{\cal M}'(0) |.\]
 This is a contradiction. 
Thus the Matsumoto torsion vanishes. By 
Lemma \ref{lemMaHo}, 
$F$ must be a Randers metric.
\qed

\bigskip
\noindent{\it Proof of Corollary \ref{cor1}}: It follows from Theorem \ref{thmC-reducible} that 
 $F$ is a Randers metric. 
Since the Matsumoto torsion vanishes, weakly Landsbergian Randers metrics are Landsberg metrics. Then the corollary follows from Numata's theorem \cite{Nu}.
\qed

\bigskip
\noindent{\it Proof of Corollary \ref{cor2}}: It follows from Theorem \ref{thmC-reducible} that $F$ is a Randers metric. According to \cite{ChSh}
a Randers metric has constant E-curvature if and only if it has constant S-curvature. 
Thus $F$ has constant S-curvature. By a theorem in \cite{Mo2}, every Finsler metric of scalar curvature with constant S-curvature must be of constant flag curvature in dimension greater than two. Thus the flag curvature of $F$ is a non-zero constant. By Numata's theorem \cite{Nu}, $F$ must be Riemannian. 
\qed

\bigskip

Clearly, Corollaries \ref{cor1} and \ref{cor2} still hold 
for noncompact Finsler manifolds of dimension $n\geq 3$ if an additional growth condition on the Matsumoto torsion as in Theorem \ref{thmC-reducible2} is imposed.

\noindent
Xiaohuan Mo\\
LMAM, School of Mathematical Sciences, Beijing University, Beijing 100871, P.R. China\\
moxh@pku.edu.cn

\bigskip

\noindent
Zhongmin Shen\\
Department of Mathematical Sciences, Indiana University-Purdue University Indianapolis, 402 N. Blackford Street, Indianapolis, IN 46202-3216, USA.  \\
zshen@math.iupui.edu

\end{document}